\RequirePackage{fix-cm}
\documentclass[smallextended]{svjour3}       
\smartqed  
\usepackage{graphicx}
\usepackage{amsfonts,amsmath,latexsym,amssymb,graphicx}
\usepackage{mathrsfs,upref}
\usepackage{mathptmx}
\begin{document}

\title{Jessen's Type Inequality and Exponential Convexity for Positive $C_0$-Semigroups
}


\author{Gul I Hina Aslam         \and
        Matloob Anwar
}


\institute{Gul I Hina Aslam \at
               School of Natural Sciences, National University of Sciences and Technology, Islamabad, Pakistan. \\
              Tel.: +92-334-8511061\\
              \email{gulihina@sns.nust.edu.pk}           
           \and
           Matloob Anwar \at
               School of Natural Sciences, National University of Sciences and Technology, Islamabad, Pakistan.
}

\date{Received: date / Accepted: date}

\maketitle

\begin{abstract}
In this paper the Jessen's type inequality for normalized positive $C_0$-semigroups is obtained.
An adjoint of Jessen's type inequality has also been derived for the corresponding adjoint-semigroup, which does not give the
analogous results but the behavior is still interesting.
Moreover, it  is followed by some results regarding positive definiteness and exponential convexity of complex structures involving operators from a semigroup.

\end{abstract}
\section{Introduction and preliminaries}
\label{intro}
A significant theory regarding inequalities and exponential convexity for real valued functions, has been developed \cite{csgji,exp}. The intention to generalize such concepts for the $C_0$-semigroup of operators, is motivated from \cite{ctmo}.\\
In the present article, we shall derive a Jessen's type inequality and the corresponding adjoint-inequality, for some $C_0$-semigroup and the adjoint-semigroup, respectively. \\
\\
The notion of Banach lattice was introduced to get a common abstract setting, within which one could talk about the ordering of elements. Therefore, the phenomena related to positivity can be generalized. It had mostly been studied in various types to spaces of real-valued functions, e.g. the space $C(K)$ of continuous functions over a compact topological space $K$, the Lebesque space $L^1(\mu)$ or even more generally the space $L^p(\mu)$ constructed over measure space $(X,\sum,\mu)$ for $1\leq p\leq\infty$. We shall use without further explanation the terms, order relation (ordering), ordered set, supremum, infimum. \\
Firstly, we shall go through the definition of vector lattice.
\begin{definition}[\cite{ops}]
Any (real) vector space $V$ with an ordering satisfying;
 \begin{description}
\item[$O_1$:] $f\leq g$ implies $f+h\leq g+h$ for all $f,g,h\in V,$
\item[$O_2$:] $f\geq 0$ implies $\lambda f\geq 0$ for al $f\in V$ and $\lambda\geq 0,$
\end{description}
is called an \emph{ordered vector space}.
\end{definition}
 The axiom $O_1$, expresses the translation invariance and therefore implies that the ordering of an ordered vector space $V$ is completely determined by the positive part $V_+=\{f\in V: f\geq 0\}$ of $V$. In other words, $f\leq g$ if and only if $f-g\in V_+$.
Moreover, the other property $O_2$, reveals that the positive part of V is a convex  set and a cone with vertex $0$ (mostly called the \emph{positive cone} of V).
\begin{itemize}
 \item An ordered vector space $V$ is called a \emph{vector lattice}, if any two elements $f,g\in V$ have a supremum, which is denoted by $sup(f,g)$ and an infimum denoted by $inf(f,g)$.\\
     It is trivially understood that the existence of supremum of any two elements in an ordered vector space implies the existence of supremum of finite number of elements in $V$. Furthermore, $f\geq g$ implies $-f\leq -g$, so the existence of finite infima therefore implied.
 \item Few important quantities are defined as follows \begin{eqnarray*}
                                   sup(f,-f) &=& |f|\quad \emph{(absolute value of $f$)} \\
                                   sup(f,0) &=& f^+ \quad \emph{(positive part of $f$)} \\
                                   sup(-f,0) &=& f^- \quad \emph{(negative part of $f$)}.
                                 \end{eqnarray*}

\item Some compatibility axiom is required, between norm and order. This is given in the following short way:
\begin{equation}\label{}
  |f|\leq |g| \quad implies\quad \|f\|\leq \|g\|.
\end{equation}
The norm defined on a vector lattice is called a lattice norm.
\end{itemize}
Now, we are in a position, to define a Banach lattice in a formal way.
\begin{definition}
 A \emph{Banach lattice} is a Banach space $V$ endowed with an ordering $\leq$, such that $(V,\leq)$ is a vector lattice with a lattice norm defined on it.\\
 A Banach lattice transforms to \emph{Banach lattice algebra}, provided $u,v\in V_+$ implies $uv\in V_+$.
 \end{definition}
 \hspace{\fill} $\square$\\
 A linear mapping $\psi$ from an ordered Banach space $V$ into itself is \emph{positive} (denoted by: $\psi\geq 0$) if $\psi f\in V_+$, for all $f\in V_+$. The set of all positive linear mappings forms a convex cone in the space $L(V)$ of all linear mappings from $V$ into itself, defining the natural ordering of $L(V)$. The absolute value of $\psi$, if it exists, is given by
$$|\psi|(f)=sup\{\psi h: |h|\leq f\},\quad (f\in V_+).$$

Thus $\psi:V\rightarrow V$ is positive if and only if $|\psi f|\leq\psi|f|$ holds for any $f\in V$.
\begin{lemma}{\textbf{[\cite{ops}, P.249]}}
A bounded linear operator $\psi$ on a Banach lattice V is a positive contraction if and only if $\|(\psi f)^+\|\leq\|f^+\|$ for all $f\in V$.
\end{lemma}
\hspace{\fill} $\square$\\
An operator $A$ on $V$ satisfies the positive minimum principle if for all $f\in D(A)_+=D(A)\cap V_+$, $\phi\in V_+'$
\begin{equation}\label{pmp}
  \langle f,\phi\rangle = 0\quad implies\quad \langle Af,\phi\rangle \geq 0.
\end{equation}
\begin{definition}
A (one parameter) $C_0$-semigroup (or strongly continuous semigroup) of operators on a Banach space $V$ is a family $\{Z(t)\}_{t\geq 0}\subset B(V)$ such that
 \begin{description}
 \item[(i)] $Z(s)Z(t)=Z(s+t)$ for all $s,t\in\mathbb{R}^+$.
 \item[(ii)] Z(0)=I, the identity operator on V.
  \item[(iii)] for each fixed $f\in V$, $Z(t)f\rightarrow f$(with respect to the norm on V) as $t\rightarrow 0^{+}$.
  \end{description}
  \end{definition}
Where $B(V)$ denotes the space of all bounded linear operators defined on a Banach space V.
\begin{definition}
The (infinitesimal) generator of $\{Z(t)\}_{t\geq 0}$ is the densely defined closed linear operator $A:V\supseteq D(A)\rightarrow R(A)\subseteq V$ such that
$$D(A)=\{f:f\in V,\lim_{t\rightarrow 0^+}A_tf\,\,exists\,in\,V\}$$
$$Af=\lim_{t\rightarrow 0^+}A_tf\,\,\,(f\in D(A))$$ where, for $t> 0$, $$A_tf=\frac{[Z(t)-I]f}{t}\,\,\,\,(f\in V).$$
\end{definition}
\hspace{\fill} $\square$\\
Let $\{Z(t)\}_{t\geq 0}$ be the strongly continuous positive semigroup, defined on a Banach lattice V. The positivity of the semigroup is equivalent to
$$|Z(t)f|\leq Z(t)|f|,\quad t\geq 0,\quad f\in V.$$
Where for positive contraction semigroups $\{Z(t)\}_{t\geq 0}$, defined on a Banach lattice V we have;
$$\|(Z(t)f)^+\|\leq \|f^+\|,\quad for\,all\,f\in V.$$
The literature presented in \cite{ops}, guarantees the existence of the strongly continuous positive semigroups and positive contraction semigroups on Banach lattice V, with some conditions imposed on the generator. The very important amongst them is, that it must always satisfy (\ref{pmp}).\\
 A Banach algebra $X$, with the multiplicative identity element $e$, is called the \emph{unital Banach algebra}. We shall call the strongly continuous semigroup $\{Z(t)\}_{t\geq 0}$ defined on $X$, a \emph{normalized semigroup}, whenever it satisfies
\begin{equation}\label{nsg}
  Z(t)(e)=e,\quad for\,all\quad t>0.
\end{equation}
The notion of normalized semigroup is inspired from normalized functionals \cite{exp}.
The theory presented in next section, is defined on such semigroups of positive linear operators defined on a Banach lattice $V$.

\section{Jessen's Type Inequality}
\label{sec:1}
In 1931, Jessen \cite{jie} gave the generalization of the Jensen's inequality for a convex function and positive linear functionals. See (\cite{buk}, PP-47). We shall prove this inequality for a normalized positive $C_0$-semigroup and convex operator, defined on a Banach lattice.\\
 Throughout the present section, $V$ will always denote a unital Banach lattice algebra, endowed with an ordering $\leq$.
\begin{definition}
Let $U$ be a nonempty open convex subset of $V$. An operator $F:U\rightarrow V$ is convex if it satisfies
\begin{equation}\label{co}
  F(tu+(1-t)v)\leq tF(u)+(1-t)F(v),
\end{equation}
whenever $u,v\in U$ and $0\leq t\leq 1$.
\end{definition}
\hspace{\fill} $\square$\\
Let $\mathfrak{D}_c(V)$ denotes the set of all differentiable convex functions $\phi:V\rightarrow V$.
\begin{theorem}{\textbf{(JESSEN'S Type INEQUALITY)}}
Let $\{Z(t)\}_{t\geq 0}$ be the positive $C_0$-semigroup on $V$ such that it satisfies (\ref{nsg}). For an operator $\phi\in\mathfrak{D}_c(V)$ and $t\geq0$;
\begin{equation}\label{jti}
  \phi(Z(t)f)\leq Z(t)(\phi f),\quad  f\in V.
\end{equation}
\end{theorem}
\noindent\textbf{Proof:}\
 Since $\phi:V\rightarrow V$ is convex and differentiable, by considering an operator-analogue of  [Theorem A, PP-98,\cite{cfs}], we have for any $f_0\in V$, there is a fixed vector $m=m(f_0)=\phi'(f_0)$ such that
 $$\phi(f)\geq \phi(f_0)+m(f-f_0),\quad f\in V.$$
 Using the property (\ref{nsg}) along with the linearity and positivity of operators in a semigroup, we obtain
\begin{equation*}
  Z(t)(\phi(f))\geq \phi(f_0)+m(Z(t)f-f_0),\quad f\in V,t\geq 0.
\end{equation*}
In this inequality, set $f_0=Z(t)f$ and the assertion (\ref{jti}) follows.
\hspace{\fill} $\square$\\

The existence of an identity element and condition (\ref{nsg}), imposed in hypothesis of the above theorem is necessary. We shall elaborate the said, by following examples.
\begin{example}
Let $X:= C_0(\mathbb{R})$, $\{Z(t)\}_{t\geq 0}$ be the left shift semigroup defined on X and $\phi$ taking the mirroring along $y$-axis. The identity function does not contain a compact support and therefore is not in $X$. If we now take a bell-shaped curve like $f(x):= e^{-x^2}$, $x\in\mathbb{R}$. Then f is positive, $\phi f=f$, and $Z(t)(\phi f)$ has maximum at $x=−t$, and it is between 0 and 1 elsewhere.
On the other hand, $\phi(Z(t)f)$ has a maximum at $s=t$ and it is immediate that we cannot compare the two functions in the usual ordering. See FIGURE 1(a).
\end{example}

\smartqed
\begin{example}
Let $\Gamma :=\{z\in\mathbb{C}:|z|=1\}$, and $X=C(\Gamma)$. The rotation semigroup $\{Z(t)\}_{t\geq 0}$ is defined as,
$T(t)f(z)=f(e^{it}\cdot z)$, $f\in X$.
The identity element $E\in X$, s.t. for all $z\in\Gamma$, $E(z)=z$. Then $Z(t)E(z)=E(e^{it}\cdot z)=e^{it}\cdot z$. Or we can say that any complex number
$z=e^{ix}$ is mapped to $e^{i(x+t)}$. $Z(t)$ satisfies (\ref{nsg}), only when $t$ is a multiple of $2\pi$. Let $f(z)=\mathfrak{R}(z)+1 >0$,
then $(Z(t)f)(e^{iz})=f(e^{i(t+z)})=cos(t+z)+1$, hence $\phi(Z(t)f)(e^{iz})=cos(t-z)+1$.
On the other hand, $\phi(f)=f$, and $Z(t)(\phi f)(e^{iz})=(Z(t)f)(e^{iz})=cos(t+z)+1$. Hence, the equality holds in (\ref{jti}) when $t$ is a multiple of $2\pi$,
but the two sides are not comparable in general. It can easily be verified that $Z'=\{Z(2\pi t)\}_{t=0}^\infty$ is a subgroup of $Z=\{Z(t)\}_{t=0}^\infty$,
as $Z(2\pi t)Z(2\pi s)=Z(2\pi (t+s))$. Therefore $Z'$ is a normalized semigroup. See FIGURE 1(b).
\end{example}
\begin{figure}[h]
 \includegraphics[width=7.1cm]{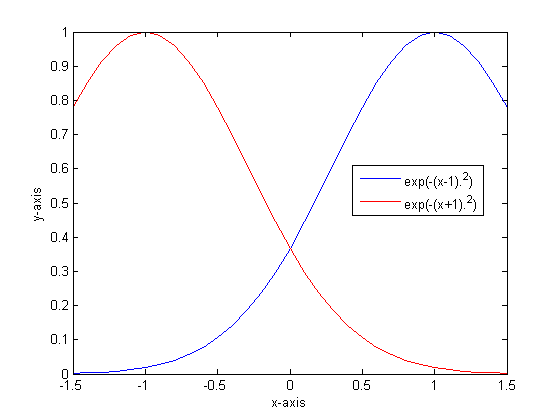}\includegraphics[width=7.1cm]{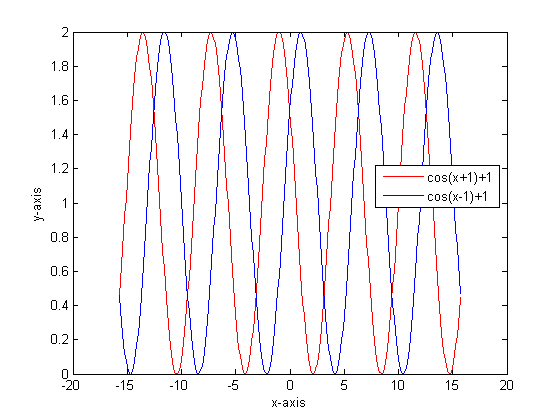}\\
\caption{(a,b)}
\end{figure}
\hspace{\fill} $\square$\\
\section{Adjoint-Jessen's Type Inequality}
In previous section, a Jessen's type inequality has been derived, for a normalized positive $C_0$-semigroup $\{Z(t)\}_{t\geq 0}$. This gives us the motivation towards knowing the behaviour of its corresponding adjoint semigroup $\{Z^\ast(t)\}_{t\geq 0}$
on $V^\ast$. As the theory for dual spaces gets more complicated, we do not expect to have the analogous results. It may ask for a detail introduction towards
a part of the dual space $V^\ast$, for which an Adjoint of Jessen's type inequality makes sense.
\begin{definition}
Given two Banach spaces $X$ and $Y$ and a bounded linear operator $L:X\rightarrow Y$, recall that the adjoint $L^\ast : Y^\ast\rightarrow X^\ast$ is defined by
\begin{equation}\label{4}
  (L^\ast y^\ast)x := y^\ast(Lx),\quad y^\ast\in Y^\ast, x\in X.
\end{equation}
\end{definition}
For a strongly continuous positive semigroup $\{Z(t)\}_{t\geq 0}$ on a Banach space $X$, by defining $Z^\ast(t)=(Z(t))^\ast$ for every $t$,
we get a corresponding adjoint semigroup $\{Z^\ast(t)\}_{t\geq 0}$ on the dual space $X^\ast$.
In \cite{aps}, it is obtained that, the adjoint semigroup $\{Z^\ast(t)\}_{t\geq 0}$ fails in general to be strongly continuous. The investigation \cite{ph}, shows that $\{Z^\ast(t)\}_{t\geq 0}$ acts in a  strongly continuous way on;
\begin{equation*}
  X^{\bigodot} := \{x^\ast\in X^\ast : Lim_{t\rightarrow 0}\|Z^\ast(t)x^\ast -x^\ast \|=0\}
\end{equation*}
This is the maximal such subspace on $X^\ast$. The space $X^{\bigodot}$ was introduced by Philips in 1955, and latter has been studied extensively
by various authors. At the present moment, we do not necessarily require the strong continuity of the adjoint semigroup $\{Z^\ast(t)\}_{t\geq 0}$ on $X^\ast$. \\
If $X$ is an ordered vector space, we say that a functional $x^\ast$ on $X$ is positive if $x^\ast (x)\geq 0$, for each $x\in X$. By the linearity of $x^\ast$,
this is equivalent to $x^\ast$ being order preserving. i.e. $x\leq y$ implies $x^\ast(x)\leq x^\ast (y)$. The set $P$ of all positive linear functionals on $X$,
is a cone in $X^\ast$.\\
We are mainly interested in the study of the space $V^\ast$, where in our case $V$ is a Banach lattice algebra.
Let us consider the regular ordering among the elements of $V^\ast$. i.e. $v_1^\ast \geq v_2^\ast$, whenever $v_1^\ast (v)\geq v_2^\ast (v)$,
for each $v\in V$. \\
Consider the convex operator (\ref{co}). In case of the equality, $F$ is simply a linear operator and the adjoint $F$ can be defined as above. But how can it be defined in other case? This question has already been answered.\\
 In \cite{pa}, some kind of adjoint has been associated to a nonlinear operator $F$. In fact, this is possible for Lipschitz continuous operators only. Consider the Banach space $\mathfrak{Lip}_0(X,Y)$ of all
Lipschitz continuous operators $F : X \rightarrow Y$ satisfying $F(\theta) = \theta$, equipped with the norm
\begin{equation*}
  [F]_{Lip}= sup_{x_1\neq x_2}\frac {\|F(x_1)-F(x_2)\|}{\|x_1-x_2\|},\quad x_1,x_2\in X.
\end{equation*}
Where $\theta\in X$ is the identity.
It is easy to see that the space $L(X, Y )$ of all bounded linear operators from X to Y is a closed subspace of $\mathfrak{Lip}_0(X,Y)$. In particular, we set
\begin{equation*}\label{}
  \mathfrak{Lip}_0(X,\mathbb{K}):= X^{\sharp}
\end{equation*}
and call $X^{\sharp}$ the pseudo-dual space of $X$; this space contains the usual dual space $X^\ast$ as closed subspace.\\
\begin{definition}
For $F\in\mathfrak{Lip}_0(X,Y)$, the pseudo-adjoint $F^{\sharp} : Y^{\sharp}\rightarrow X^{\sharp}$ of $F$ is defined by;
\begin{equation}\label{}
  F^{\sharp}(y^\sharp)(x) := y^\sharp (F(x)),\quad y^\sharp \in Y^\sharp, x\in X.
\end{equation}
\end{definition}
This is of course a straightforward generalization of (\ref{co}); in fact, for linear operators $L$ we have $L^\sharp|_{Y^\ast} = L^\ast$. i.e. the restriction of the
pseudo-adjoint to the dual space is the classical adjoint. \\
For the sake of convenience, we shall denote the adjoint of the operator $F$ by $F^\ast$, throughout the present section. Either it's a classical adjoint or the Pseudo-adjoint (depending upon the operator $F$).\\
Similarly, the considered dual space of the vector lattice algebra $V$ will be denoted by $V^\ast$, which can be the intersection of the pseudo-dual and classical dual spaces in case of a nonlinear convex operator.
\begin{lemma}
Let $F$ be the convex operator on a Banach space $X$, then the adjoint operator $F^\ast$ on the dual space $X^\ast$ is also convex.
\end{lemma}
\noindent\textbf{Proof:}\
For $x\in X$ and $0\leq \lambda\leq 1$\begin{eqnarray*}
               (F^\ast(\lambda x_1^\ast +(1-\lambda)x_2^\ast),x) &=& (\lambda x_1^\ast +(1-\lambda)x_2^\ast , F(x)),\quad \\
                &=& \lambda(x_1^\ast , F(x))+(1-\lambda)(x_2^\ast,F(x)),
             \end{eqnarray*}
where $x_1^\ast,x_2^\ast\in X^\ast$. By putting $x=\mu x+(1-\mu)x$, for $0\leq \mu\leq 1$ and using the convexity of the operator $F$ we finally get
\begin{equation*}
  F^\ast(\lambda x_1^\ast +(1-\lambda)x_2^\ast)\leq \lambda F^\ast(x_1^\ast)+(1-\lambda)F^\ast(x_2^\ast).
\end{equation*}
Hence, $F^\ast$ is convex on $X^\ast$.

\begin{theorem}{\textbf{(ADJOINT-JESSEN'S INEQUALITY)}}
Let $\{Z^\ast(t)\}_{t\geq 0}$ be the adjoint semigroup on $V^\ast$ such that the original semigroup $\{Z(t)\}_{t\geq 0}$, the operator $\phi$ and the space $V$ are same
as in Theorem (1). For a convex operator $\phi^\ast: V^\ast\rightarrow V^\ast$ and $t\geq 0$
\begin{equation}\label{ajti}
  \phi^\ast(Z^\ast(t)f^\ast)\geq Z^\ast(t)(\phi^\ast f^\ast),\quad  f^\ast\in V^\ast.
\end{equation}
\end{theorem}
\noindent\textbf{Proof:}\
  For $f\in V$ and $t\geq 0$, consider
  \begin{eqnarray*}
    (\phi^\ast [Z^\ast (t)f^\ast],f) &=& (Z^\ast(t)f^\ast,\phi(f)) \\
     &=& (f^\ast,Z(t)(\phi f)) \\
     & \geq & (f^\ast,\phi(Z(t)f)) \\
     &=& (\phi^\ast (f^\ast),Z(t)f) \\
     &=& (Z^\ast(t)[\phi^\ast f^\ast],f)
  \end{eqnarray*}
  Therefore, the assertion (\ref{ajti}) is satisfied.

\section{Exponential Convexity}
In this section we shall define the exponential convexity of an operator. Moreover, few complex structures, involving the operators from a semigroup, will be proved to be exponentially convex.
\begin{definition}
Let $V$ be a Banach lattice endowed with ordering $\leq$. An operator $H:I\rightarrow V$ is exponentially convex if it is continuous and for all $n\in\mathbb{N}$
\begin{equation}\label{}
  \sum_{i,j=1}^n\xi_i\xi_j H(x_i+x_j)f\geq 0,\quad f\in V,
\end{equation}
where $\xi_i\in\mathbb{R}$ such that  $x_i+x_j\in I\subseteq\mathbb{R}$, $1\leq i,j\leq n$.
\end{definition}
\begin{proposition}
Let $V$ be a Banach lattice endowed with ordering $\leq$. For an operator $H:I\rightarrow V$, the following propositions are equivalent.
\begin{description}
   \item[(i)] $H$ is exponentially convex.
    \item[(ii)] $H$ is continuous and for all $n\in\mathbb{N}$
\begin{equation}\label{}
  \sum_{i,j=1}^n\xi_i\xi_j H\big(\frac{x_i+x_j}{2}\big)f\geq 0,\quad f\in V,
\end{equation}
where $\xi_i\in\mathbb{R}$ and $x_i\in I\subseteq\mathbb{R}$, $1\leq i\leq n$.
 \end{description}
\end{proposition}
\noindent\textbf{Proof:}\
$(i)\Rightarrow (ii)$\\
Take any $\xi_i\in\mathbb{R}$ and $x_i\in I$, $1\leq i\leq n$. Since the interval $I\subseteq\mathbb{R}$ is convex, the midpoints, $\frac{x_i+x_j}{2}\in I$. Now set $y_i=\frac{x_i}{2}$, for $1\leq i\leq n$. Then we have, $y_i+y_j=\frac{x_i+x_j}{2}\in I$, for all $1\leq i,j\leq n$. Therefore, for all $n\in\mathbb{N}$, we can apply $(i)$ to get,
\begin{equation*}
  \sum_{i,j=1}^n\xi_i\xi_j H(y_i+y_j)f=\sum_{i,j=1}^n\xi_i\xi_j H\big(\frac{x_i+x_j}{2}\big)f\geq 0,\quad f\in V.
\end{equation*}
$(ii)\Rightarrow (i)$\\
Let $\xi_i,x_i\in\mathbb{R}$, such that $x_i+x_j\in I$, for $1\leq i,j\leq n$. Define $y_i=2x_i$, so that $x_i+x_j=\frac{y_i+y_j}{2}\in I$. Therefore, for all $n\in\mathbb{N}$, we can apply $(ii)$ to get,
\begin{equation*}
  \sum_{i,j=1}^n\xi_i\xi_j H\big(\frac{y_i+y_j}{2}\big)f = \sum_{i,j=1}^n\xi_i\xi_j H(x_i+x_j)f\geq 0,\quad f\in V.
\end{equation*}
\hspace{\fill} $\square$\\

\begin{remark}
Let $H$ be an exponentially convex operator. Writing down the fact for $n=1$, in (10), we get that $H(x)f\geq 0$, for $x\in I$ and $f\in V$. For $n=2$, we have
\begin{equation*}
  \xi_1^2H(x_1)f+2\xi_1\xi_2 H\big(\frac{x_1+x_2}{2}\big)f+\xi_2^2H(x_2)f\geq 0.
\end{equation*}
Hence, for $\xi_1=-1$ and $\xi_2=1$, we have
\begin{equation*}
  H\big(\frac{x_1+x_2}{2}\big)f\leq \frac{H(x_1)f+H(x_2)f}{2},
\end{equation*}
i.e. $H:I\rightarrow V$, does indeed satisfy the condition of convexity.
\end{remark}
\hspace{\fill} $\square$\\
For $U\subseteq V$, let us assume that $F:U\rightarrow V$ is continuously differentiable on $U$. i.e. the mapping $F': U\rightarrow \mathfrak{L}(V)$, is continuous. Moreover $F''(f)$, will be a continuous linear transformation from $V$ to $\mathfrak{L}(V)$. A bilinear transformation $B$ defined on $V\times V$ is symmetric if $B(f,g)=B(g,f)$ for all $f,g\in V$. Such a transformation is \textbf{positive definite [nonnegative definite]}, if for every nonzero $f\in V$, $B(f,f)>0$ [$B(f,f)\geq 0$]. Then, $F''(f)$ is symmetric wherever it exists. See [\cite{cfs}, pp-69].
\begin{theorem}[[\cite{cfs}, P.100]]\label{condif}
Let $F$ be continuously differentiable and suppose that second derivative exists throughout an open convex set $U\subseteq V$. Then $F$ is convex on $U$ if and only if $F''(f)$ is nonnegative definite for each $f\in U$. And if $F''(f)$ is positive definite on $U$, then $F$ is strictly convex.
\end{theorem}
\hspace{\fill} $\square$\\
\begin{definition}[\cite{log}]
 Let $V$ be a Banach algebra with unit $e$. For $f\in V$, we define a function $log(f)$ from $V$ to $V$.
\begin{equation*}
  log(f)= -\sum_{n=1}^\infty\frac{(e-f)^n}{n}= -(e-f)-\frac{(e-f)^2}{2}-\frac{(e-f)^3}{3}-...
\end{equation*}
for $||(e-x)||\leq 1$.
\end{definition}
\begin{lemma}
 Let $V$ be a unital Banach algebra. For $f\in V$, a family of operators $F_t$ is defined as
\begin{equation}
 F_t(f)=\begin{cases}
                 \frac{f^t}{t(t-1)}, & t\neq 0,1; \\

                 -\log f, & t=0;\\

                  f\log f, & t=1.
               \end{cases}
\end{equation}
Then $D^2F_t(f):= f^{t-2}$. Whenever, $f\in V_+$, $D^2F_t(f)\in V_+$, therefore by Theorem 3, the mapping $f\rightarrow F_t(f)$ is convex.
\end{lemma}
\begin{theorem}\label{t4}
Let $\{Z(t)\}_{t\geq 0}$ be the positive $C_0$-semigroup, defined on a unital Banach lattice algebra $V$, such that it satisfies (\ref{nsg}). Let $f\in V$, such that $f^r\in V$, for $r\in \mathbb{R}\supseteq I\ \{0,1\}$, $log f\in V$, if $r=0$ and $f log f\in V$, if $r=1$. Let us define
\begin{equation}
  \Lambda_t := Z(t)(F_t(f))-F_t(Z(t)f)
\end{equation}
Then\begin{description}
      \item[(i)] for every $n\in\mathbb{N}$ and for every $p_k\in I$, $k=1,2,...,n$,
\begin{equation}\label{a}
  \big[\Lambda_{\frac{p_i+p_j}{2}}\big]_{i,j=1}^n\geq 0.
\end{equation}
      \item[(ii)] If the mapping $f\rightarrow\Lambda_t$ is continuous on $I$, then it is exponentially convex on $I$.
    \end{description}
\end{theorem}
\noindent\textbf{Proof:}\
Consider the operator
\begin{equation*}
  G(f)=\sum_{i,j=1}^n u_iu_j F_{p_{ij}}(f)
\end{equation*}
for $f>0$, $u_i\in\mathbb{R}$ and $p_{ij}\in I$ where $p_{ij}=\frac{p_i+p_j}{2}$. Then
\begin{equation*}
  D^2G(f):= \sum_{i,j=1}^nu_iu_jf^{p_{ij}-2}=\big(\sum_{i=1}^1 u_if^{\frac{p_i}{2}-1}\big)^2\geq 0,\quad f>0.
\end{equation*}
So, $G(f)$ is a convex operator. Therefore by applying (\ref{jti}) we get
\begin{equation*}
  \sum_{i,j=1}^n u_iu_j\Lambda_{p_{ij}}\geq 0,
\end{equation*}
and the assertion (\ref{a}) follows. Assuming the continuity and using the \textbf{Proposition 1} we have also exponential convexity of the operator $f\rightarrow\Lambda_t$.
\hspace{\fill} $\square$\\
\begin{lemma}
Let $V$ be a unital Banach algebra, for $f\in V$, let us define the following family of operators
\begin{equation*}
  H_t(f)=\begin{cases}
                 \frac{e^{tf}}{t^2}, & t\neq 0; \\
\\
                 \frac{f^2}{2}, & t=0.
               \end{cases}
\end{equation*}
Then, $D^2H_t(f)=e^{tf}$. By Theorem (\ref{condif}), the mapping $f\rightarrow H_t(f)$, is convex on $V$.
\end{lemma}
\begin{theorem}
For $\Lambda_t := Z(t)(H_t(f))-H_t(Z(t)f)$, $(i)$ and $(ii)$ from Theorem (\ref{t4}), holds.
\end{theorem}
\noindent\textbf{Proof:}\ Similar to the proof of Theorem (\ref{t4}).
\hspace{\fill} $\square$\\
\\
\textbf{Competing Interests}\\
The authors declare that they have no competing interests.\\

\noindent\textbf{Author's Contribution}\\
All authors contributed equally and significantly in writing this paper. All authors read and approved the final manuscript.\\

\noindent\textbf{Acknowledgement}\\
 Authors of this paper are grateful to Prof. Andr\'{a}s B\'{a}tkai for his generous help in construction of examples.

\end{document}